%% file: membrane_locking_regge.tex
\documentclass[10pt,twocolumn,fleqn]{article}
\usepackage[numbers]{natbib}

\usepackage{ifthen}
\usepackage{stmaryrd}
\usepackage{amscd}
\usepackage{bm}

\usepackage{ifthen}
\usepackage{stmaryrd}
\usepackage{amscd}
\usepackage{bm}
\usepackage{amssymb}
\usepackage{amsmath}
\usepackage{color}
\usepackage{cleveref}
\usepackage{pgfplots}

\usepackage{chngcntr}
\counterwithin{equation}{section}
\counterwithin{figure}{section}
\counterwithin{table}{section}

\newcommand{\R}{\ensuremath{\mathbb{R}}}
\newcommand{\mat}[1]{\ensuremath{\begin{pmatrix}#1\end{pmatrix}}}
\newcommand{\T}{\ensuremath{\mathcal{T}}}
\newcommand{\E}{\ensuremath{\mathcal{E}}}
\newcommand{\tr}[1]{\ensuremath{\,\text{tr}(#1)}}
\newcommand{\sym}[1]{\ensuremath{\,\text{sym}(#1)}}
\newcommand{\re}[1]{\ensuremath{\hat{#1}}}
\newcommand{\HCurl}[1][]{\ensuremath{H(\text{curl}\ifthenelse{\equal{#1}{}}{}{,#1})}}
\newcommand{\Hone}[1][]{\ensuremath{H^1\ifthenelse{\equal{#1}{}}{}{(#1)}}}
\newcommand{\Cinf}[1][]{\ensuremath{C^{\infty}\ifthenelse{\equal{#1}{}}{}{(#1)}}}
\newcommand{\Ltwo}[1][]{\ensuremath{L^2\ifthenelse{\equal{#1}{}}{}{(#1)}}}
\newcommand{\idop}{\text{id}}
\newcommand{\bsigma}{\bm{\sigma}}
\newcommand{\bF}{\bm{F}}
\newcommand{\bM}{\bm{M}}
\newcommand{\bI}{\bm{I}}
\newcommand{\bP}{\bm{P}}
\newcommand{\bq}{\bm{q}}
\newcommand{\RegInt}[1][]{\mathcal{I}^{R}_{h\ifthenelse{\equal{#1}{}}{}{,#1}}}
\newcommand{\HoneInt}[1][]{\mathcal{I}_{h\ifthenelse{\equal{#1}{}}{}{,#1}}}
\newcommand{\nv}{\nu}
\newcommand{\Proj}{\bP_{\tau}}
\newcommand{\W}{\ensuremath{\mathcal{W}}}
\newcommand{\CGS}{\bm{C}_{\tau}}
\newcommand{\GS}{\bm{E}_{\tau}}
\newcommand{\Regge}{\text{Reg}}

\begin{document}
\begin{center}
{\Large \textbf{Avoiding Membrane Locking with Regge Interpolation}} \\[0.5cm]
\end{center}
Michael Neunteufel$^{a,}$\footnote{Corresponding author.\\E-mail address: michael.neunteufel@tuwien.ac.at }, Joachim Sch{\"o}berl$^{a}$\\
$^a$ Institute for Analysis and Scientific Computing, TU Wien, Wiedner Hauptstra{\ss}e 8-10, 1040 Wien, Austria

\begin{abstract}
In this paper we present a novel method to overcome membrane locking of thin shells. An interpolation operator into the so-called Regge finite element space is inserted in the membrane energy term to weaken the implicitly given kernel constraints. The number of constraints is asymptotically halved on triangular meshes compared to reduced integration techniques. Provided the interpolant, this approach can be incorporated easily into any shell element. The performance of the proposed method is demonstrated by means of several benchmark examples.\newline

\end{abstract}

\section{Introduction}
\label{sec:introduction}
In mathematical formulations of plates and shells a small parameter, the thickness $t$, is involved. The lack of finite element approximations fulfilling the implicitly  given constraints of the physical model leads to so-called locking phenomena \cite{babuvska92,Chenais1994}. As the thickness becomes small --depending heavily on the geometry and boundary conditions-- the shell falls in one of two different categories: the membrane dominated or bending dominated case \cite{CHAPELLE98}. For shells, shear and membrane locking can be observed in the case of non-inhibited pure bending. The former, induced by the Kirchhoff constraint in the limit case, has been extensively discussed and a variety of shear locking free plate and shell elements have been proposed. For membrane locking, also called inextensional locking, where the (curved) elements fail to represent pure bending, only little numerical analysis has been done \cite{arnold97,koschnick05,gerdes98,choi98,CS98}. A framework based on discrete models for constructing low order membrane locking free elements has been proposed in \cite{Quag16}. In practice, mostly reduced integration schemes combined with stabilization techniques \cite{zienk71,stolarski82,stolarski83,PS86} and assumed strain methods \cite{McNeal82,PS86,HH86,JP87,Chap11} are used. Therein, the membrane constraints are weakened due to under-integration, without deteriorating the membrane stability in the membrane dominated case, or the strain components are evaluated at certain points and extrapolated, respectively. The discrete strain gap method \cite{koschnick05}, related to assumed strain procedures, modifies the normal strains eliminating parasitic membrane parts. Further, mixed methods introducing the membrane force tensor as additional unknown have been proposed \cite{arnold97,CS98,EOB13}. It is well known that p and hp-refinement strategies \cite{pitkaranta92,suri96,HLP96} may overcome the problem of membrane locking, but in the case of low order triangular elements only little influence of the reduced integration techniques has been observed \cite{choi98}. 

Tullio Regge derived in \cite{Regge61} a geometric discretization of the Einstein field equations by approximations with a piece-wise constant metric. In theoretical, and later also numerical, physics so-called Regge calculus was applied e.g., in fields of relativity and quantum mechanics  and has been further developed the last fifty years, see \cite{williams92} for an overview. An analytical perspective of Regge calculus was given in \cite{cheg86,Cheeger84}. It has been observed that Regge's approach is equivalent to specify lengths at all edges of a mesh, analogical to Whitney-forms \cite{whitney57}. In the context of finite element exterior calculus (FEEC) \cite{arnold06, arnold10} a finite element point of view was given in \cite{christiansen04,christiansen11} and the resulting Regge elements have been generalized to arbitrary polynomial order on triangles and tetrahedrons \cite{li18}.

In this work the resulting Regge interpolant is used to construct membrane locking free methods for shells, staying stable in the case of inhibited pure bending. It can be incorporated into any existing method and finite element code, provided the interpolation operator. In a variety of numerical examples the performance of the method is tested.

This paper is structured as follows: In the next section we will give an overview of the construction of Regge elements and the corresponding interpolation operator. Section 3 is devoted to the description of the proposed method in the general setting of shells. In Section 4 the method is discussed and in the last section we apply the resulting algorithm to several established membrane locking benchmark examples.

\section{Regge elements}
\label{sec:regge_elements}
For the reader's convenience we give a brief introduction in the construction of Regge finite elements. Then the implementation of the Regge interpolation operator is discussed. 

As the Regge elements approximate symmetric tensor fields, we seek for a matrix valued finite element space. To prescribe the edge lengths globally only the tangential-tangential components need to be continuous. Therefore, let $\T_h$ a triangulation of a bounded domain $\Omega\subset\R^2$ and $\E_h$ the corresponding skeleton, i.e., the set of all edges. The set of all piece-wise polynomials up to degree $k$ on $\T_h$ and $\E_h$ is denoted by $\Pi^k(\T_h)$ and $\Pi^k(\E_h)$, respectively. For each element $T\in\T_h$ the tangential and outer normal vector on the boundary $\partial T$ are given by $t$ and $n$, respectively. The outer dyadic product of two vectors is denoted by $\otimes$.

The Regge finite element space is given by
\begin{flalign}
& \Regge^k:=\{\bsigma\in[\Pi^k(\T_h)]^{2\times 2}_{\text{sym}}:\llbracket\bsigma_{tt}\rrbracket=0  \},&&
\end{flalign}
where $\bsigma_{tt}:= (\bI-n\otimes n)\bsigma(\bI-n\otimes n)$ is the tangential-tangential component of $\bsigma$, $\bI$ the identity matrix, and $\llbracket\bsigma_{tt}\rrbracket$ denotes the corresponding jump over elements.

Further, we define the Lagrangian nodal finite element space as
\begin{flalign}
& V^k_h:=\Pi^k(\T_h)\cap C^0(\Omega),&&
\end{flalign}
where $C^0(\Omega)$ denotes to set of all continuous functions on $\Omega$. For a (high-order) construction of $\Hone$-conforming finite elements we refer to \cite{Br2013,ZT20,Zaglmayr06}.

In the context of \cite{hiptmair99} we define functionals, the degrees of freedom (dofs), on the reference or physical triangle $T$ with the local space $[\Pi^k(T)]^{2\times 2}_{\text{sym}}$. Therefore, let $\{q_{E,l}\}$ and $\{\bq_{T,l}\}$ denote a polynomial basis of $\Pi^k(E_{ij})$ on the edge $E_{ij}$ between two vertices, $i\neq j$, and $[\Pi^{k-1}(T)]^{2\times 2}_{\text{sym}}$ on the triangle $T$, respectively. Then the functionals read
\begin{flalign}
\Psi_{E_{ij},l}:&\,\bsigma\mapsto \int_{E_{ij}}\bsigma:q_{E,l}t_E\otimes t_E\,ds,&&\label{eq:functional_edge}\\
\Psi_{T,l}:&\,\bsigma\mapsto \int_{T}\bsigma:\bq_{T,l}\,dx,\label{eq:functional_vol}&&
\end{flalign}
where $\bm{A}:\bm{B}:=\sum_{ij}\bm{A}_{ij}\bm{B}_{ij}$ denotes the Frobenius scalar product and $t_E$ the tangent vector of the edge $E_{ij}$. Note that $t_E\otimes t_E$ is single valued, i.e., does not depend on the orientation of $t_E$.\\

\begin{figure}
\centering
\includegraphics[width=0.2\textwidth]{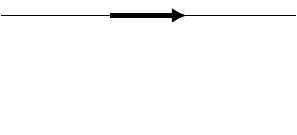}
\includegraphics[width=0.2\textwidth]{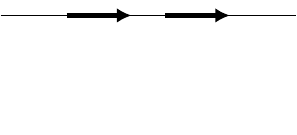}

\includegraphics[width=0.2\textwidth]{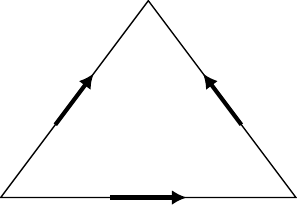}
\includegraphics[width=0.2\textwidth]{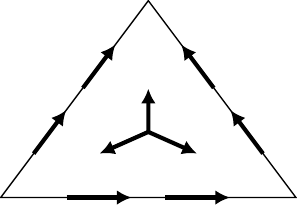}

\caption{Lowest order and first order Regge elements for segments and triangles.}
\label{fig:regge_el_segm_trig}
\end{figure}

In one dimension the Regge elements coincide with $\Ltwo$-conforming discontinuous finite elements. On triangular elements the dofs are associated with the edges \eqref{eq:functional_edge}, analogical to $\HCurl[]$-conforming \cite{Nedelec1980,Nedelec1986} elements, and inner bubbles \eqref{eq:functional_vol} for higher polynomial degrees. The lowest order (polynomial order $k=0$) and first order Regge elements on segments and triangles are illustrated in Figure \ref{fig:regge_el_segm_trig}.

We now give an explicit basis for the corresponding shape functions $\{\bm{\varphi}_i\}$ to \eqref{eq:functional_edge}--\eqref{eq:functional_vol}. Therefore, let $\lambda_{i}\in\Pi^1(\hat{T})$ denote the barycentric coordinates of the vertices $V_1=(-1,0)$, $V_2=(1,0)$, and $V_3=(0,1)$ of the reference triangle $\hat{T}$, i.e. $\lambda_{i}(V_j)=\delta_{ij}$, where $\delta_{ij}$ denotes the Kronecker delta. Then, the lowest order basis functions are given by
\begin{flalign}
\bm{\varphi}_{E_{ij},0}:= \nabla\lambda_{i}\odot\nabla\lambda_{j},&&\label{eq:basis_edges}
\end{flalign}
where $\odot$ denotes the symmetric dyadic product of two vectors. For a generalization to arbitrary order $k$ let for $\alpha,\beta > -1$, $p_n^{(\alpha,\beta)}$ and $\hat{p}_n^{(\alpha,\beta)}$ denote the $n$-th Jacobi and $n$-th integrated Jacobi polynomial \cite{A74,AAR99,BS06}, respectively, 
\begin{flalign}
p_n^{(\alpha,\beta)}(x)&:=\frac{1}{2^nn!(1-x)^{\alpha}(1+x)^{\beta}}&&\nonumber\\ &\times\frac{d^n}{dx^n}((1-x)^{\alpha}(1+x)^{\beta}(x^2-1)^n),\quad n\in\mathbb{N}_0,&&\\
\hat{p}_n^{(\alpha,\beta)}(x)&:=\int_{-1}^xp_{n-1}^{(\alpha,\beta)}(\zeta)\,d\zeta,\qquad n\geq 1,\,\hat{p}_0^{(\alpha,\beta)}(x)=1.&&
\end{flalign}
As we are only interested in the case $\beta=0$ the notation can be simplified by $p_n^{(\alpha,0)}(x)=p_n^{\alpha}(x)$ and $\hat{p}_n^{(\alpha,0)}(x)=\hat{p}_n^{\alpha}(x)$. The latter fulfill the following orthogonality properties 
\begin{flalign}
&\int_{-1}^1(1-x)^{\alpha}p_j^{\alpha}(x)p_l^{\alpha}(x)\,dx = \frac{2^{\alpha+1}}{2j+\alpha+1}\delta_{jl},&&\\
&\int_{-1}^1(1-x)^{\alpha}\hat{p}_j^{\alpha}(x)\hat{p}_l^{\alpha}(x)\,dx = 0\qquad \text{for }|j-l| > 2.&&
\end{flalign}
Note that with $\alpha=0$ the (integrated) Jacobi polynomials reduce to the (integrated) Legendre polynomials. Hence, the high order edge basis functions read, $l=1,\dots,k$,
\begin{flalign}
\label{eq:basis_edges_ho}
\bm{\varphi}_{E_{ij},l}:= \hat{p}^0_l\left(\frac{\lambda_{i}-\lambda_{j}}{\lambda_{i}+\lambda_{j}}\right)(\lambda_{i}+\lambda_{j})^l\nabla\lambda_{i}\odot\nabla\lambda_{j}.&&
\end{flalign}
The (high order) cell basis functions of order $k>0$ are given by, $l_1,l_2 \geq 0$,
\begin{subequations}
\label{eq:basis_inner}
\begin{flalign}
&\bm{\varphi}_{T^1,l_1,l_2}:=w^{(l_1,l_2)}\,\lambda_{1}(\nabla\lambda_{2}\odot\nabla\lambda_{3}),\quad l_1+l_2 \leq k-1,&&\\
&\bm{\varphi}_{T^2,l_1,l_2}:=w^{(l_1,l_2)}\,\lambda_{2}(\nabla\lambda_{3}\odot\nabla\lambda_{1}),\quad l_1+l_2 \leq k-1,&&\\
&\bm{\varphi}_{T^3,l_1,l_2}:=w^{(l_1,l_2)}\,\lambda_{3}(\nabla\lambda_{1}\odot\nabla\lambda_{2}),\quad l_1+l_2 \leq k-1,&&
\end{flalign}
\end{subequations}
with the Dubiner basis
\begin{flalign}
w^{(l_1,l_2)}:=p^0_{l_1}\left(\frac{x}{1-y}\right)(1-y)^{l_1}p_{l_2}^{2l_1+1}(2y-1).&&
\end{flalign}
The shape functions \eqref{eq:basis_edges}, \eqref{eq:basis_edges_ho}, and \eqref{eq:basis_inner} build a basis of $[\Pi^k(\T_h)]^{2\times 2}_{\text{sym}}$, which has a dimension of $3(k+1)(k+2)/2$: The edge shape functions are linearly independent as $(\bm{\varphi}_{E_{i},l})_{t_{E_j}t_{E_j}}=\delta_{ij}$, where $E_i$ and $t_{E_i}$ denotes the i-th edge and corresponding tangent vector. Further, $(\bm{\varphi}_{T^i,l_1,l_2})_{t_{E_j}t_{E_j}}=0$ for $i,j=1,2,3$ and thus, the inner shapes are independent of the edge basis. The claim follows together with the independence of $\lambda_{i}(\nabla\lambda_{j}\odot\nabla\lambda_{k})$ for $i\neq j\neq k$ and counting all shape functions. A different basis for triangular Regge elements is constructed in \cite{li18}.

Given functionals \eqref{eq:functional_edge}--\eqref{eq:functional_vol} $\{\Psi_i\}$ and the corresponding shape functions $\{\bm{\varphi}_i\}$, one can define the following Regge interpolation operator
\begin{flalign}
\RegInt[k]&:[\Cinf[\Omega]]^{2\times 2}_{\text{sym}}\rightarrow \Regge^k,&&\nonumber\\
&\bsigma\mapsto\sum_{i=0}^{N_k}\alpha_i\bm{\varphi}_i,&&
\end{flalign}
where $N_k\in\mathbb{N}$ denotes the number of degrees of freedom (the number of shape functions), $\Cinf[\Omega]$ the set of all smooth functions on $\Omega$, and the coefficients $\alpha_i$ are obtained by the following consideration:

\begin{figure}
	\centering
	\begin{tabular}{cc}
	\includegraphics[width=0.22\textwidth]{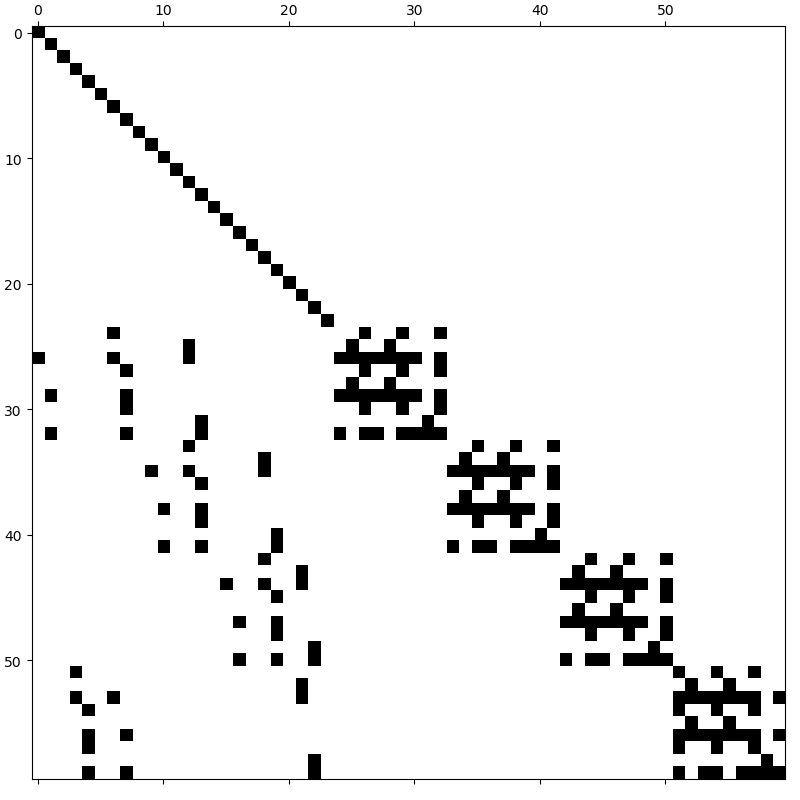}&
	\includegraphics[width=0.22\textwidth]{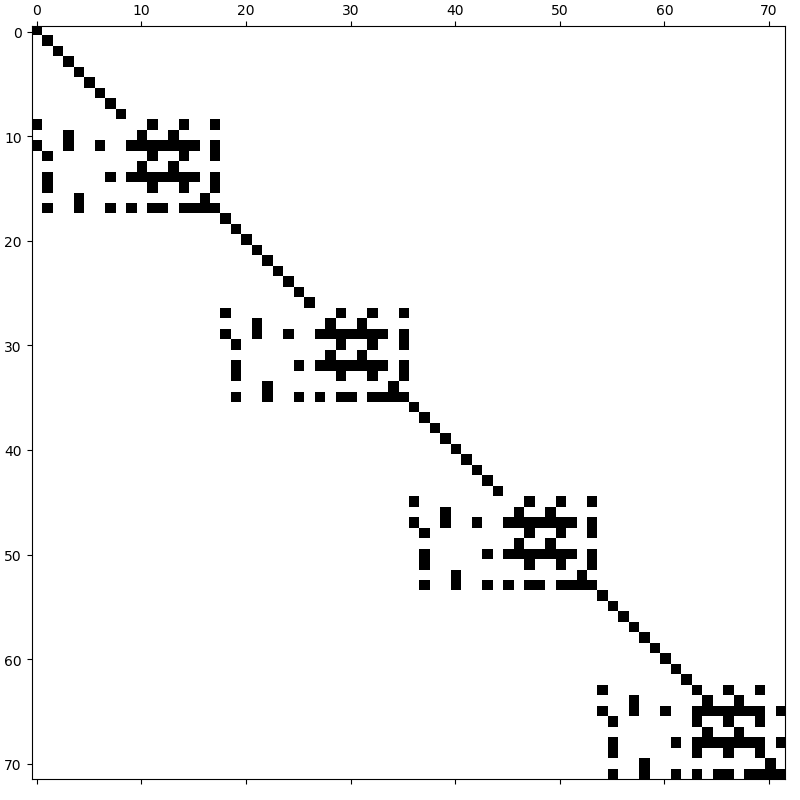}\\
	(a) & (b)
	\end{tabular}
	
	\caption{Sparsity pattern of dual mass matrix $\bM^{D}$ on the unit square divided by its diagonals, i.e., consisting of 4 triangles. (a) Sparsity pattern of Regge elements of order $k=2$. (b) Sparsity pattern of Regge elements of order $k=2$ with broken continuity.}
	\label{fig:sparsity_pattern}
\end{figure}

Let $\bsigma_h=\RegInt[k](\bsigma)$ be the interpolation of a given function $\bsigma$. Then, with \eqref{eq:functional_edge}--\eqref{eq:functional_vol}, $\bsigma_h$ is the solution of the following system of equations
\begin{flalign}
& \sum_{E\in\E_h}\int_E\bsigma_h:q_{\E}t_E\otimes t_E\,ds = \sum_{E\in\E_h}\int_E\bsigma:q_{\E}t_E\otimes t_E\,ds, &&\nonumber\\
& \sum_{T\in\T_h}\int_T\bsigma_h:\bq_{\T}\,dx = \sum_{T\in\T_h}\int_T\bsigma:\bq_{\T}\,dx, &&\label{eq:equation_dualmass}
\end{flalign}
for all $\bq_{\T}\in[\Pi^{k-1}(\T_h)]^{2\times 2}_{\text{sym}}$ and $q_{\E}\in\Pi^{k}(\E_h)$. In matrix form we obtain the linear equation
\begin{flalign}
\bM^{D} \langle\alpha_i\rangle = f,\,\, \bM^{D} =\mat{\bM^{D}_{EE} & \bM^{D}_{ET}\\ \bM^{D}_{TE} &\bM^{D}_{TT}},\,\, f = \mat{f_E\\ f_T},&&
\end{flalign}
with $\langle\alpha_i\rangle$ denoting the corresponding coefficient vector of $\bsigma_h$ and $\bM^{D}_{EE} = \Psi_{E_{ij}}(\bm{\varphi}_{E_{ij}})$, $\bM^{D}_{ET} = \Psi_{E_{ij}}(\bm{\varphi}_{T})$, $f_E = \Psi_{E_{ij}}(\bsigma)$ and analogously for the other components.
There holds by construction $\Psi_{E_{ij}}(\bm{\varphi}_{T})=0$ and thus, $\bM^{D}_{ET}=0$, resulting in a lower block triangular dual mass matrix, see Figure \ref{fig:sparsity_pattern} (a). This structure can be exploited for inverting $\bM^{D}$. Note that the matrix is not symmetric, as in general $\Psi_{T}(\bm{\varphi}_{E_{ij}})\neq0$. As will be discussed later the tangential-tangential continuity of the elements can be broken leading to a block diagonal matrix depicted in Figure \ref{fig:sparsity_pattern} (b). Note that in this case the dofs are ordered element-wise instead of splitting them into edge and inner dofs.\\

Given a Regge element $\hat{\bsigma}_h$ on the reference triangle $\hat{T}$, one needs to map it to a physical element $T=\Phi(\hat{T})$, $\Phi:\hat{T}\rightarrow\R^2$, in such a way that the tangential-tangential continuity is preserved. For $\HCurl$-conforming elements the covariant transformation $u\circ\Phi:=\bF^{-\top}\hat{u}$, $\bF :=\nabla \Phi$, is used to preserve the tangential continuity. Thus, by applying the covariant transformation on both sides
\begin{flalign}
\bsigma_h\circ\Phi:=\bF^{-\top} \hat{\bsigma}_h\bF^{-1}\label{eq:trafo_hcurlcurl_vol}&&
\end{flalign}
leads to a symmetric tangential-tangential continuous function $\bsigma_h$ on $T$. For shells one needs to map the reference triangle $\hat{T}$ to a (possibly curved) surface triangle $T\subset\R^3$. The deformation gradient $\bF\in\Cinf(\Omega,\R^{3\times 2})$ is not invertible and thus, the transformation rule \eqref{eq:trafo_hcurlcurl_vol} needs to be adopted by using the Moore--Penrose pseudo-inverse instead
\begin{flalign}
\bsigma_h\circ\Phi:=(\bF^{\dagger})^\top \hat{\bsigma}_h\bF^{\dagger}\label{eq:trafo_hcurlcurl_surf}.&&
\end{flalign}

It is possible to obtain a geometry free dual mass matrix $\bM^D$ by transforming the basis $q_{\E}$ and $\bq_{\T}$ in \eqref{eq:equation_dualmass} with
\begin{flalign}
&q_{\E}\circ\Phi:=J_b\hat{q}_{\E}, &&\\
&\bq_{\T}\circ\Phi:=\frac{1}{J}\bF\hat{\bq}_{\T}\bF^\top,&&\label{eq:trafo_vol_dual}
\end{flalign}
where $J_b:= \|\bF t\|_2$ denotes the boundary determinant and $J:=\det(\bF)$, leading to ($T=\Phi(\hat{T})$, $E=\Phi(\hat{E})$)
\begin{flalign}
\int_E\bsigma_h:q_{\E}t\otimes t\,ds &= \int_{\hat{E}}J_b\bF^{-\top}\hat{\bsigma}_h\bF^{-1}:\frac{\hat{q}_{\E}}{J_b}(\bF \hat{t})\otimes(\bF \hat{t})\,d\hat{s}&& \nonumber\\
&= \int_{\hat{E}}\hat{\bsigma}_h:\hat{q}_{\E}\hat{t}\otimes \hat{t}\,d\hat{s},&&\\
\int_T\bsigma_h:\bq_{\T}\,dx &=\int_{\hat{T}}J\bF^{-\top}\hat{\bsigma}_h\bF^{-1}:\frac{1}{J}\bF\hat{\bq}_{\T}\bF^\top\,d\hat{x}&&\nonumber\\
&=\int_{\hat{T}}\hat{\bsigma}_h:\hat{\bq}_{\T}\,d\hat{x}.&&
\end{flalign}
Therefore, one may exploit this property to assemble only (permutations of) one element and ``fill'' the whole matrix $\bM^D$. This procedure works also for transformation \eqref{eq:trafo_hcurlcurl_surf}, where one has to choose the surface determinant $J=\sqrt{\det(\bF^\top\bF)}$ for $\bF\in \R^{3\times 2}$ or $J=\|\text{cof}(\bF)\nu\|_2$ for the extended matrix $\bF\in \R^{3\times 3}$ in \eqref{eq:trafo_vol_dual} with $\text{cof}(\bF)$ denoting the cofactor matrix and $\nu$ the surface normal vector.

\section{Methodology}
\subsection{Shells and membrane energy}
Let $\re{\Omega}\subset\R^3$ be an undeformed configuration of a shell with thickness $t$, described by its mid-surface $S$ and the corresponding orientated normal vector $\nv$
\begin{flalign}
\Omega:=\{\re{x}+z\nv(\re{x}): \re{x}\in S, z\in [-t/2,t/2]\}.&&
\end{flalign}
Furthermore, let $\Phi:\Omega\rightarrow\R^3$ be the deformation from the initial to the deformed configuration of the shell and $\phi:S\rightarrow\R^3$ the deformation of the approximated mid-surface. The corresponding triangulation consisting of possibly curved triangles of $S$ is denoted by $\T_h$. Then, we define $\bF:=\nabla_{\tau}\phi$ as the deformation gradient. Here, $\nabla_{\tau} \phi$ denotes the surface gradient of $\phi$, which can be introduced in weak sense \cite{de13}, or directly as Fr\'{e}chet-derivative. We can split the deformation into the identity function and the displacement, $\phi=\idop+u$, and thus, $\bF= \Proj+\nabla_{\tau} u$ with the projection onto the tangent plane $\Proj:=\bI-\nv\otimes\nv$.\\

The shell energy functional can be split into a membrane, bending and shear energy part, cf. \cite{Ciarlet05,Chap11,BBWR04},
\begin{flalign}
\W(u) = \frac{t}{2}\,E_{\text{mem}}(u)+\frac{t^3}{2}E_{\text{bend}}(u)+\frac{t}{2}E_{\text{shear}}(u)-f(u),&&\label{eq:shell_prob}
\end{flalign}
where $f$ denotes the external forces.

We focus on the membrane energy and consider the full nonlinear term 
\begin{flalign}
\label{eq:mem_end_nonlin}
&E_{\text{mem}}(u):=\int_{S}\|\GS\|^2_{\bM}\,dx,&&
\end{flalign}
with $\GS:=1/2(\CGS-\bI)$ denoting the Green strain tensor restricted on the tangent plane, $\CGS:= \bF^\top\bF$ the Cauchy--Green strain tensor. The material norm is given by 
\begin{flalign}
\|\cdot\|^2_{\bM}:=\frac{\bar{E}}{1-\bar{\nu}^2}\left(\bar{\nu}\tr{\cdot}^2+(1-\bar{\nu})\tr{\cdot^2}\right),&&
\end{flalign}
with the material tensor $\bM$, the Young's modulus $\bar{E}$, and the Poisson's ratio $\bar{\nu}$, respectively.

The linearization of \eqref{eq:mem_end_nonlin} is given by
\begin{flalign}
\label{eq:mem_end_lin}
& E^{\text{lin}}_{\text{mem}}(u):=\int_{S}\|\sym{\Proj\nabla_{\tau}u}\|^2_{\bM}\,dx,&&\nonumber\\
& (E^{\text{lin}}_{\text{mem}}(u))_{\alpha\beta}=\int_{S}\|\frac{1}{2}(u_{\alpha|\beta}+u_{\beta|\alpha})-b_{\alpha\beta}u_3\|^2_{\bM}\,dx,&&
\end{flalign}
where $u_{\alpha|\beta}$ denotes the covariant derivative, $b_{\alpha\beta}$ the second fundamental form, $\alpha,\beta \in\{1,2\}$, and $u_3$ the displacement component in normal direction, see e.g. \cite{Chap11} for the notation of curvilinear coordinates.

\subsection{Usage of Regge interpolant}
In what follows let the discrete displacements $u_h \in [V_h^k]^3$. For the proposed method we insert the Regge interpolation operator of order $k-1$ into the membrane energy \eqref{eq:mem_end_nonlin}
\begin{flalign}
&\int_{S}\|\RegInt[k-1] \GS\|^2_{\bM}\,dx.\label{eq:reg_inter_in_memen}&&
\end{flalign}
Due to the tangential-tangential continuity of $\GS$ --the discrete Jacobian $\nabla_{\tau} u_h$ is tangential-continuous-- it is mathematically equivalent to apply the projection operator $\RegInt[k-1]$ only element-wise. This makes the method cheap, as no additional global system has to be solved. Further, properties as symmetry and positivity gets preserved. Therefore, only small problems of the form \eqref{eq:equation_dualmass} on one element have to be solved on each integration point. 

Further, \eqref{eq:reg_inter_in_memen} is equivalent to a three-field formulation by introducing the (locally and thus discontinuous) Regge interpolant $\bm{R}\in \Regge^{k-1,\mathrm{dc}}$ and the corresponding local shape functionals $\bm{Q}\in[\Regge^{k-1,\mathrm{dc}}]^*$ as additional unknowns. Note that the functionals \eqref{eq:functional_edge}--\eqref{eq:functional_vol} span the topological dual space $[\Regge^{k-1,\mathrm{dc}}]^*$. The corresponding Lagrangian reads
\begin{flalign}
&\mathcal{L}(u,\bm{R},\bm{Q}):=\int_{S}\|\bm{R}\|^2_{\bM}\,dx + \langle \bm{R} - \GS,\bm{Q}\rangle_{\T_h},\label{eq:memlock_regge_threefield}&&
\end{flalign}
where, according to \eqref{eq:functional_edge}--\eqref{eq:functional_vol},
\begin{flalign}
&\langle \bm{R}-\GS,\bm{Q}\rangle_{\T_h}:=\sum_{T\in\T_h}\Big(\int_T(\bm{R}- \GS):\bm{Q}_{T}\,dx\nonumber\\
&\qquad\quad+\sum_{E\in\partial T}\int_{E}(\bm{R}- \GS):Q_{E}t_{E}\otimes t_{E}\,ds\Big)&&
\end{flalign}
and thus, with \eqref{eq:equation_dualmass}, $\bm{R}=\RegInt[k-1]\GS$. In numerical experiments we observed that if the full nonlinear Green strain tensor $\GS$ is used in \eqref{eq:memlock_regge_threefield} less Newton iterations are needed than for the direct interpolation procedure \eqref{eq:reg_inter_in_memen}.

\section{Discussion}
Let $u$ be the exact solution of the shell problem \eqref{eq:shell_prob} in the case of non-inhibited pure bending such that $E^{\text{lin}}_{\text{mem}}(u)=0$. Interpolating $u$ into the Lagrangian finite element space $[V_h^k]^3$, $u_h := \HoneInt[k] u$, $\HoneInt$ denoting the standard nodal interpolation operator, does not guarantee in general that $E^{\text{lin}}_{\text{mem}}(u_h)=0$ for the discrete displacements. I.e., the interpolation operator does not preserve the kernel of the membrane operator. Therefore, pure bending modes induce discrete membrane energy modes due to the discrete constraints. This effect dominates for small thickness parameters $t$, the shell element is called to be too stiff and locking occurs.

By using the Regge interpolant $\RegInt[k-1]E^{\text{lin}}_{\text{mem}}(u_h)$ we weak the discrete constraints. Reduced integration schemes follow the same idea, using less Gau\ss-integration points, which corresponds to an $\Ltwo$ instead of a Regge interpolation. When we compare the number of dofs, which can be interpreted as the number of constraints, one can observe that on a single triangle $T$ the number of constraints are equal, as the dimension of both spaces are the same, $\text{dim} =3(k+1)(k+2)/2$.

For a triangulation $\T_h$, however,  the number of constraints differ already in the lowest order case significantly. For Regge elements we have one degree of freedom per edge, whereas in the reduced integration scheme one has three per element. Asymptotically there holds 
\begin{flalign}
\#T \approx 2\#V,\qquad \#E \approx 3\#V,&&
\end{flalign}
where $\#T$, $\#E$, and $\#V$ denote the number of triangles, edges, and vertices of the triangulation $\T_h$, respectively. Therefore,
\begin{flalign}
\#E \approx 3\#V < 6\#V \approx 3\#T &&
\end{flalign}
and thus, the Regge interpolation reduces the number of constraints asymptotically by a factor of two compared to the $\Ltwo$-projection. Furthermore, on a triangulation $\T_h$ of a flat two-dimensional domain or a surface described by one single embedding (and thus not closed) there holds
\begin{flalign}
3+\#E=3\#V-\#V_B=2\#V+\#V_I,&&
\end{flalign}
where $\#V_B$ and $\#V_I$ denote the number of vertices on the boundary and in the inner domain of the surface, respectively.
The discrepancy of three corresponds to the number of rigid-body motions in two dimensions, two translations and one rotation. Therefore, for given displacements at the vertices one can find a unique value per edge describing the (tangential-tangential) stretching between two vertices. This fits perfectly to the following (linear) exact sequence

\begin{flalign}
  \begin{CD}
     RB @> \text{id} >> [\Cinf[\Omega]]^2 @> \nabla_{\text{sym}} >> [\Cinf[\Omega]]^{2\times 2}_{\text{sym}} \\
     && @V \HoneInt[k] VV   @V \RegInt[k-1] VV\\
     RB @> \text{id} >> [V_h^k]^2 @> \nabla_{\text{sym}} >> \Regge^{k-1}
  \end{CD},&&
\end{flalign}
where $RB:=\{Ax+b\,|\, A\in \R^{2\times 2},\, A^\top=-A, \, b\in\R^2\}$ denotes the set of linearized rigid body motions.

In \cite{christiansen11,Hauret13} they used this sequence in three dimensions as a part of a larger complex and proofed in the lowest order case commuting and exactness properties. For a nonlinear complex one has to replace the symmetric gradient by the Green strain tensor and $RB=\{Ax+b\,|\, A\in SO(2), \, b\in\R^2\}$, where $SO(2)$ denotes the set of all orthogonal $2\times 2$ matrices with determinant one.

In case of the full nonlinear membrane energy term \eqref{eq:mem_end_nonlin} the Green strain operator $\GS:[\Pi^k(\T_h)]^d\rightarrow[\Pi^{2k-2}(\T_h)]^{d\times d}_{\text{sym}}$ doubles the polynomial degree asymptotically element-wise, with the exception $k=1$. This may lead to even worse discrete kernel conservation. Thanks to the Regge interpolant, the Green strain tensor gets projected back to polynomial degree $k-1$ and again the number of constraints are reduced.

The idea of inserting an interpolation operator has already been successfully applied to avoid shear locking. E.g. for the mixed interpolated tensorial components (MITC) elements \cite{bathe87,brezzi91} an $\HCurl[]$ interpolant is inserted into the shear energy term. Also methods where the rotations get directly approximated by $\HCurl[]$-conforming finite elements overcome shear locking \cite{PS17}.\newline

In the lowest order case $k=1$ for the displacements, membrane locking is not observed as long as an isoperimetric mapping for the shell geometry is considered. Curving the geometry by a higher polynomial degree as the displacements leads to enormous membrane locking in the lowest order case. However, using the Regge interpolation $\RegInt[0]$ reduces this locking phenomena too.

\section{Numerical examples}
\label{sec:numerical examples}
To avoid shear locking effects we use the Kirchhoff--Love shell model introduced in \cite{NS19}, where the Regge interpolation has been successfully used for triangular meshes. The method is implemented in the open source finite element library Netgen/NGSolve\footnote{www.ngsolve.org} \cite{Sch97,Sch14}. 

For the benchmarks we use second order finite elements for the displacements, where the geometry is mapped isoperimetrically, i.e. curved elements are used, called method p2.

The forces are chosen such that the deformations are in the linear regime. Therefore, the differences between the linearized \eqref{eq:mem_end_lin} and full nonlinear \eqref{eq:mem_end_nonlin} membrane energy is marginal. Further, the forces are scaled appropriately with the thickness parameter $t$ ($t^3$ in the bending dominated and $t$ in the membrane dominated case) such that the deformations are in the same magnitude. Due to the nonlinear membrane and bending energy, however, the results may vary little with respect to the thickness parameter. The reference values are computed by using fourth order finite elements for the displacement on the finest mesh, called method p4, and the relative error is computed by \texttt{|result - reference|/|reference|}.

\subsection{Cylinder with free ends}
\label{subsec:ne_cyl}
\begin{figure}
\centering
\includegraphics[width=0.25\textwidth]{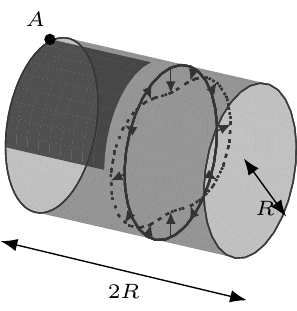}
\caption{Geometry for Cylinder with free ends benchmark.}
\label{fig:geom_cyl_shell}
\end{figure}

\begin{figure}
\centering
\includegraphics[width=0.238\textwidth]{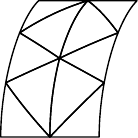}
\includegraphics[width=0.238\textwidth]{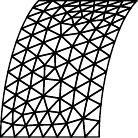}
\caption{Unstructured meshes with 10 and 160 elements for cylinder with free ends benchmark.}
\label{fig:mesh_start_cyl_shell}
\end{figure}

\begin{figure}
	\centering
	\includegraphics[width=0.238\textwidth]{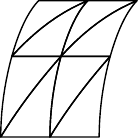}
	\includegraphics[width=0.238\textwidth]{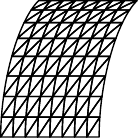}
	\caption{Structured meshes with 8 and 128 elements for cylinder with free ends benchmark.}
	\label{fig:mesh_start_cyl_shell_struct}
\end{figure}

A cylinder with free ends is loaded with a periodic force \cite{PITKARANTA95,Chap11}, see Figure \ref{fig:geom_cyl_shell}. By symmetries the computational domain is one eighth of the original and symmetry boundary conditions are used, see Figure \ref{fig:mesh_start_cyl_shell} and \ref{fig:mesh_start_cyl_shell_struct}. The material and geometric parameters are $R=1$, $E=3\times 10^4$, $\nu=0.3$, $t\in\{0.1,0.01,0.001,0.0001\}$ and the cylinder is loaded by the normal pressure distribution $P=t^3\cos(2\zeta)\re{\nu}$, $\zeta$ and $\re{\nu}$ denoting the circumferential arc-length and the normal vector on the reference configuration, respectively, cf. Figure \ref{fig:geom_cyl_shell}.


\begin{figure}
	\centering
	\includegraphics[width=0.238\textwidth]{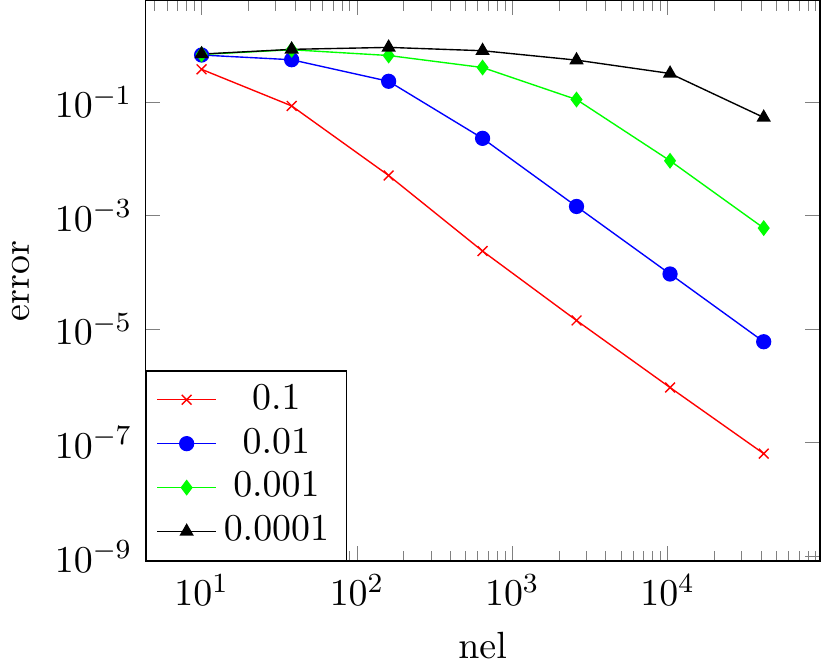}
	\includegraphics[width=0.238\textwidth]{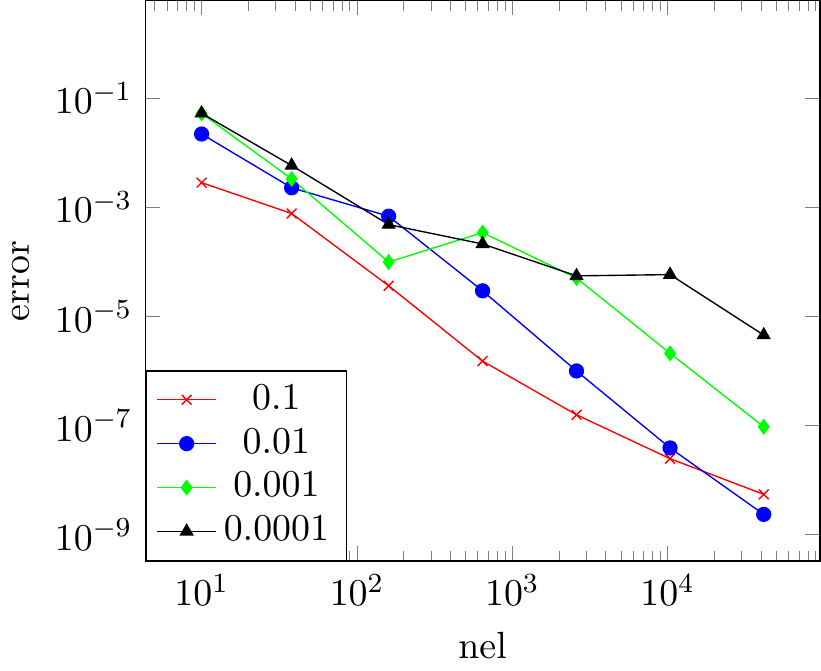}
	\caption{Results for cylinder with free ends, method p2 (unstructured mesh) without and with Regge interpolation.}
	\label{fig:results_cyl_p2_lin}
\end{figure}

\begin{figure}
	\centering
	\includegraphics[width=0.238\textwidth]{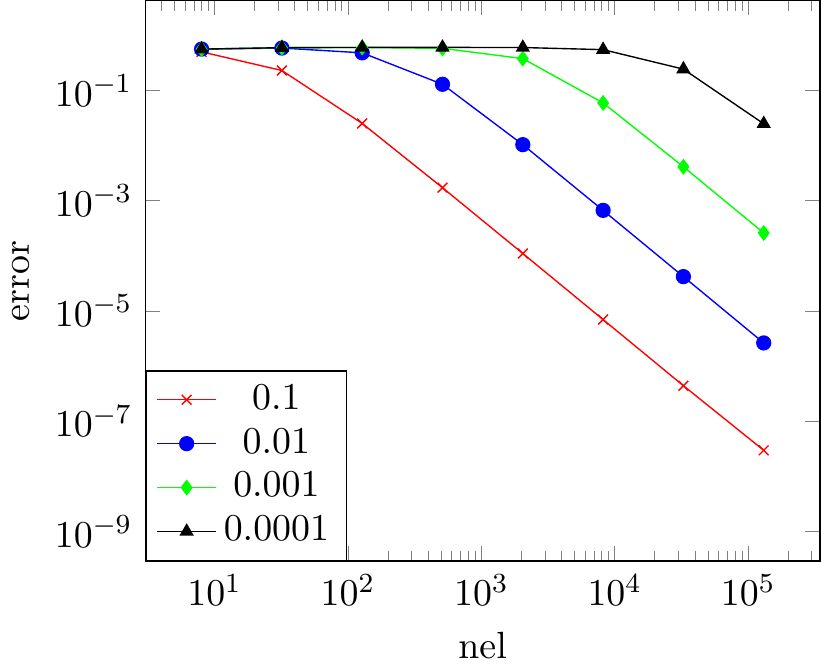}
	\includegraphics[width=0.238\textwidth]{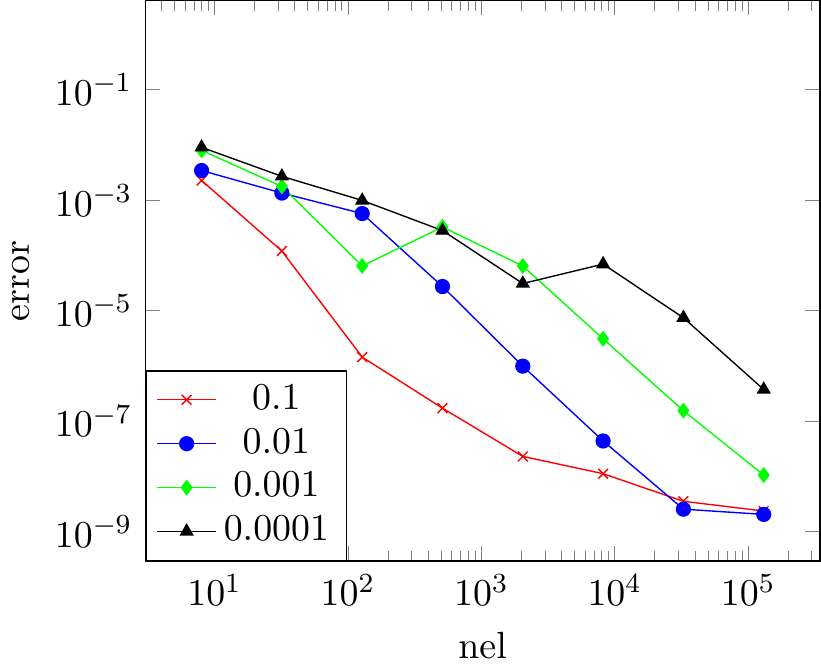}
	\caption{Results for cylinder with free ends, method p2 (structured mesh) without and with Regge interpolation.}
	\label{fig:results_cyl_p2_lin_struct}
\end{figure}

The radial deflection at point $A$ is measured and listed in Tables \ref{tab:results_cyl_unstruct_p2_nregge}-\ref{tab:results_cyl_struct_p2_regge}. The relative error for unstructured meshes can be found in Figure \ref{fig:results_cyl_p2_lin} and for structured meshes in Figure \ref{fig:results_cyl_p2_lin_struct}. There the classical locking behavior can be observed if the Regge interpolant is not used, as the pre-asymptotic range increases rapidly for smaller thicknesses. Using Regge interpolation avoids this pre-asymptotic behavior. Further, for a small amount of elements the relative errors start already with $<10$ percent, also for thick parameters $t$. The results on the unstructured meshes are comparable to the structured one.

\begin{table}
	\input{CylShell_p1_strcFalse_linTrue_regFalse.tab}
	\caption{Results for cylinder with free ends $\times 10^5$, method p2 (unstructured mesh) without Regge interpolation.}
	\label{tab:results_cyl_unstruct_p2_nregge}
\end{table}
\begin{table}
	\input{CylShell_p1_strcFalse_linTrue_regTrue.tab}
	\caption{Results for cylinder with free ends $\times 10^5$, method p2 (unstructured mesh) with Regge interpolation.}
	\label{tab:results_cyl_unstruct_p2_regge}
\end{table}
\begin{table}
	\input{CylShell_p1_strcTrue_linTrue_regFalse.tab}
	\caption{Results for cylinder with free ends $\times 10^5$, method p2 (structured mesh) without Regge interpolation.}
	\label{tab:results_cyl_struct_p2_nregge}
\end{table}
\begin{table}
	\input{CylShell_p1_strcTrue_linTrue_regTrue.tab}
	\caption{Results for cylinder with free ends $\times 10^5$, method p2 (structured mesh) with Regge interpolation.}
	\label{tab:results_cyl_struct_p2_regge}
\end{table}

\subsection{Axisymmetric hyperboloid with free ends}
\label{subsec:ne_axis_hyp}
\begin{figure}
\centering
\includegraphics[width=0.25\textwidth]{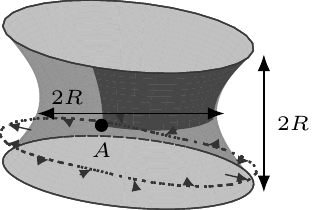}	
\caption{Geometry for axisymmetric hyperboloid with free ends benchmark.}
\label{fig:geom_hyp_shell}
\end{figure}

\begin{figure}
	\centering
	\includegraphics[width=0.238\textwidth]{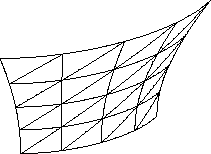}
	\includegraphics[width=0.238\textwidth]{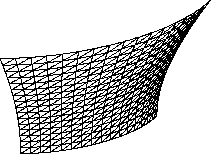}
	\caption{Meshes with 32 and 512 elements for axisymmetric hyperboloid with free ends benchmark.}
	\label{fig:mesh_start_hyp_shell}
\end{figure}

An axisymmetric hyperboloid is described by the equation 
\begin{flalign}
x^2+y^2=R^2+z^2,\quad z\in [-R,R]
\end{flalign}
with free boundaries is loaded by a force, see \cite{Chap11}. Due to symmetries it is sufficient to use one eighth of the geometry and symmetry boundary conditions, see Figure \ref{fig:geom_hyp_shell} for the geometry and Figure \ref{fig:mesh_start_hyp_shell} for a coarse and fine mesh. The material and geometric parameters are $R=1$, $E=2.85\times 10^4$, $\nu=0.3$, $t\in\{0.1,0.01,0.001,0.0001\}$, $P=\frac{t^3}{\sqrt{x^2+y^2}}\cos(2\zeta)\mat{x\\y\\0}$ similar to the previous benchmark.

\begin{figure}
	\centering
	\includegraphics[width=0.238\textwidth]{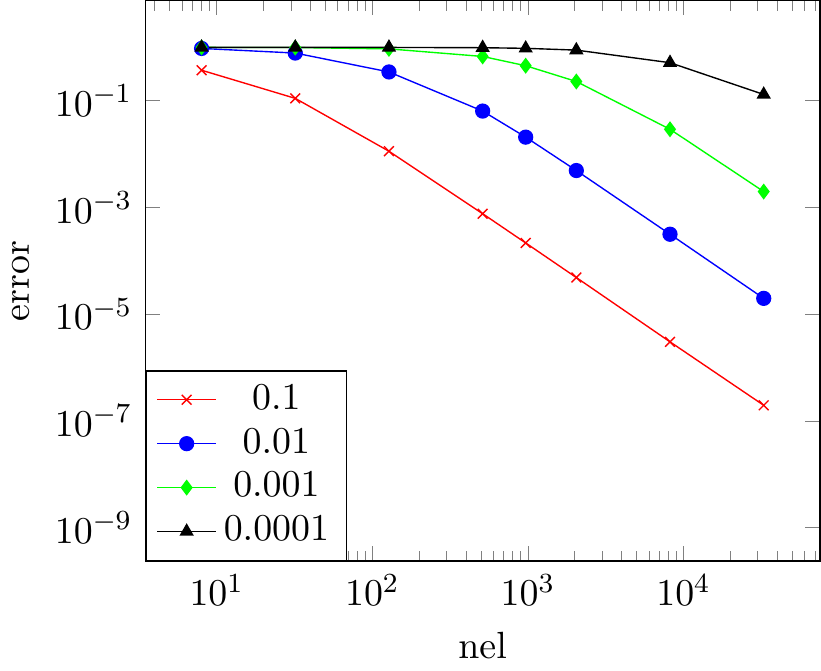}
	\includegraphics[width=0.238\textwidth]{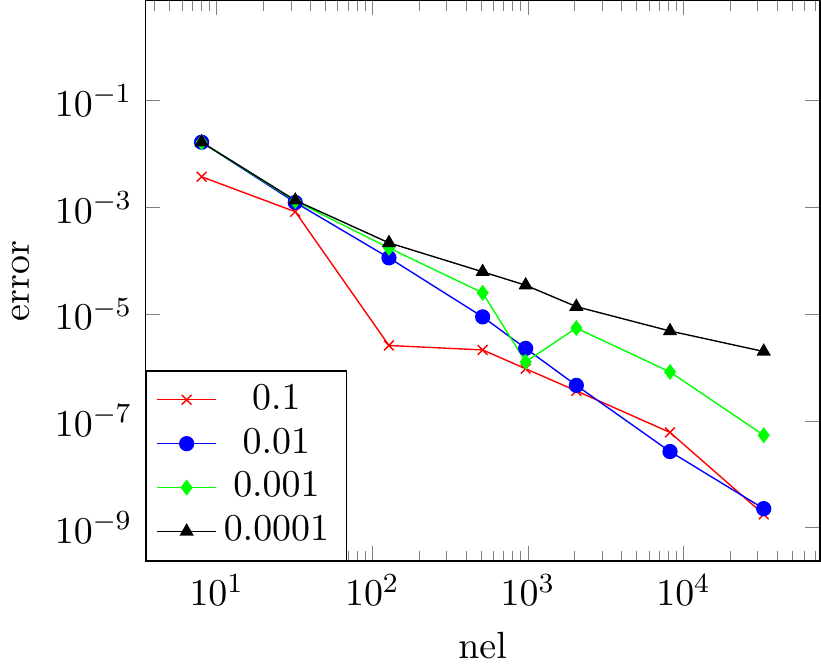}
	\caption{Results for axisymmetric hyperboloid with free ends, method p2 without and with Regge interpolation.}
	\label{fig:results_hyp_p2_lin}
\end{figure}

The radial deflection at point $A$ is listed in Table  \ref{tab:results_hyp_p2_nregge}-\ref{tab:results_hyp_p2_regge} and the relative error in Figure \ref{fig:results_hyp_p2_lin}. Again the results improve using Regge interpolation and we emphasize that for $t=0.0001$ with 8 elements the difference with a factor of $10^5$ is immensely ($-2\times 10^{-10}$ vs $-2\times 10^{-5}$ with the reference value $-1.89\times 10^{-5}$).

\begin{table}
	\input{Hyperboloid_p1_strcTrue_linTrue_regFalse.tab}
	\caption{Results for axisymmetric hyperboloid with free ends $\times 10^5$, method p2 without Regge interpolation.}
	\label{tab:results_hyp_p2_nregge}
\end{table}
\begin{table}
	\input{Hyperboloid_p1_strcTrue_linTrue_regTrue.tab}
	\caption{Results for axisymmetric hyperboloid with free ends $\times 10^5$, method p2 with Regge interpolation.}
	\label{tab:results_hyp_p2_regge}
\end{table}

\subsection{Uniform bending of cylindrical shell}
\label{subsec:ne_uni_bend_cyl}
\begin{figure}
	\centering
	\includegraphics[width=0.15\textwidth]{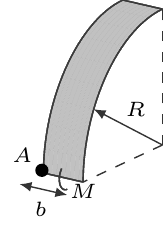}    
	\caption{Geometry for uniform bending of cylindrical shell benchmark.}
	\label{fig:geom_uni_bend_cyl_shell}
\end{figure}

\begin{figure}
	\centering
	\hspace*{0.1cm}\includegraphics[width=0.238\textwidth]{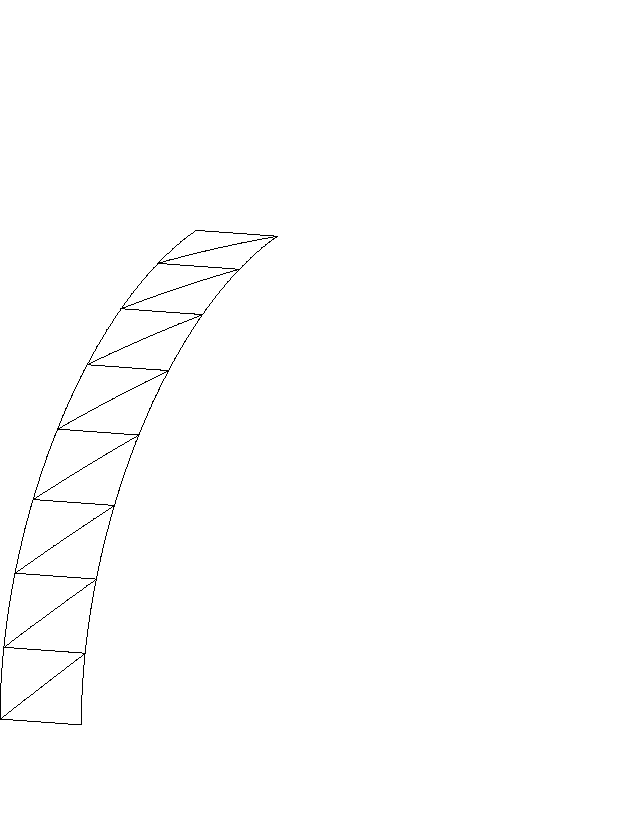}
	\includegraphics[width=0.238\textwidth]{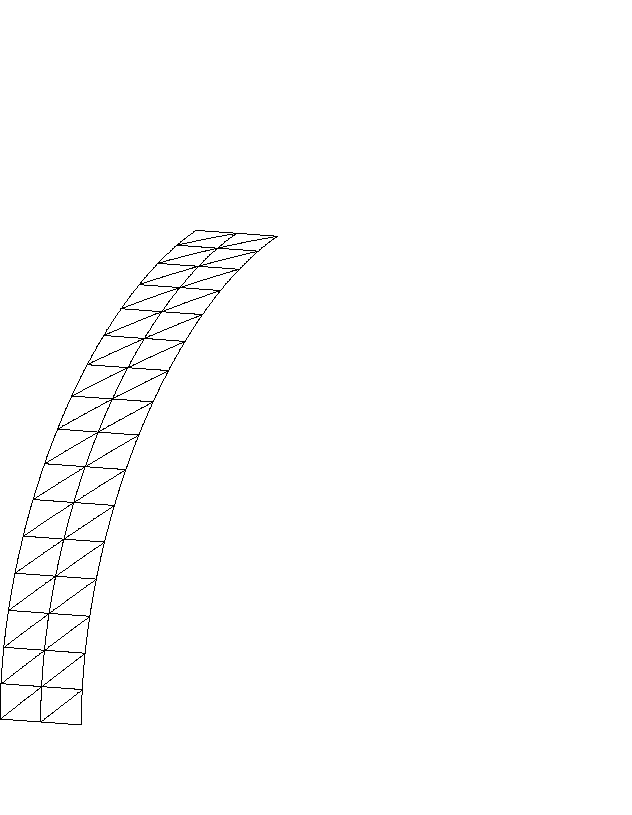}
	\caption{Meshes with 16 and 64 elements for uniform bending of cylindrical shell benchmark.}
	\label{fig:mesh_start_unibend_shell}
\end{figure}

A moment $M$ is applied to a cylindrical shell, which is fixed at the top \cite{koschnick05}. The material and geometric parameters are $R=0.1$, $b=0.025$, $E=2\times 10^5$, $\nu=0$, $t\in\{0.1,0.01,0.001,0.0001\}$, $M_0=(t/R)^3$, see Figure \ref{fig:geom_uni_bend_cyl_shell} and \ref{fig:mesh_start_unibend_shell}.

\begin{figure}
	\centering
	\includegraphics[width=0.238\textwidth]{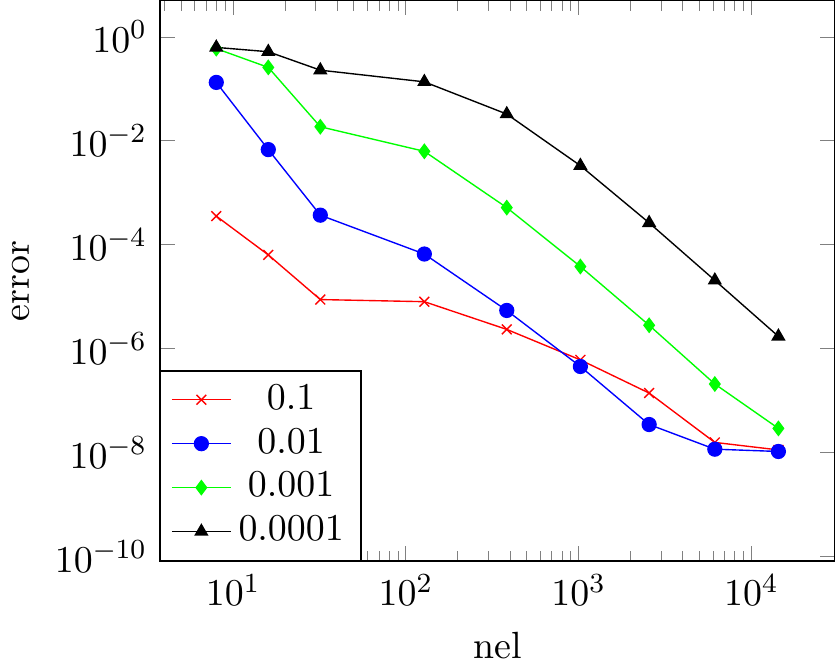}
	\includegraphics[width=0.238\textwidth]{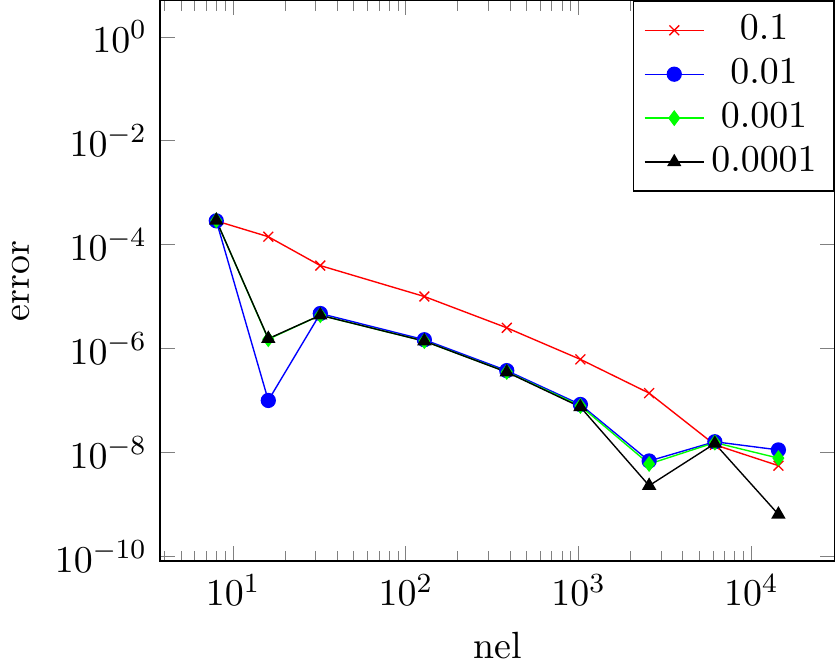}
	\caption{Results for uniform bending of cylindrical shell, method p2 without and with Regge interpolation.}
	\label{fig:results_unibend_p2_lin}
\end{figure}

This time the deflection orthogonal to the radial direction is computed at point $A$. The results can be found in Table \ref{tab:results_unibend_p2_nregge}-\ref{tab:results_unibend_p2_regge} and Figure \ref{fig:results_unibend_p2_lin}. In this benchmark the method without interpolation operator does not produce a strong pre-asymptotic regime for small thicknesses. However, the initial relative error gets larger. In contrast, the errors with the Regge interpolation start all at nearly the same value and show a uniform convergence behavior.

\begin{table}
	\input{UnifBendCyl_p1_strcTrue_linTrue_regFalse.tab}
	\caption{Results for uniform bending of cylindrical shell $\times 10^4$, method p2 without Regge interpolation.}
	\label{tab:results_unibend_p2_nregge}
\end{table}
\begin{table}
	\input{UnifBendCyl_p1_strcTrue_linTrue_regTrue.tab}
	\caption{Results for uniform bending of cylindrical shell $\times 10^4$, method p2 with Regge interpolation.}
	\label{tab:results_unibend_p2_regge}
\end{table}

\subsection{Hyperbolic paraboloid}
\label{subsec:ne_hyp_par}
\begin{figure}
\centering
\includegraphics[width=0.2\textwidth]{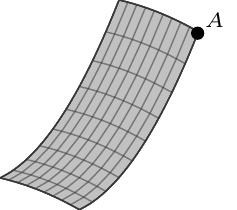}	
\caption{Geometry for hyperbolic paraboloid benchmark.}
\label{fig:geom_hyp_par_shell}
\end{figure}

\begin{figure}
	\centering
	\includegraphics[width=0.238\textwidth]{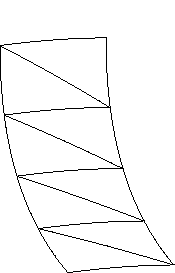}
	\includegraphics[width=0.238\textwidth]{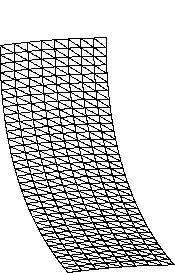}
	\caption{Meshes with 8 and 512 elements for hyperbolic paraboloid benchmark.}
	\label{fig:mesh_start_hyp_par_shell}
\end{figure}

A hyperbolic paraboloid, which is described by the embedding
\begin{flalign}
\Phi&:[0,3]\times[0,1]\rightarrow\R^3&&\nonumber\\
& (x,y)\mapsto (x,y,\alpha(y^2-x^2)),
\end{flalign}
is clamped at the bottom and subjected to a surface force $f$ \cite{choi98}. On the right side symmetry boundary conditions are used, the other boundaries are free. The material and geometric parameters are $\alpha=0.2$, $E=2.85\times 10^4$, $\nu=0.3$, $t\in\{0.1,0.01,0.001,0.0001\}$, $f=8t^3\re{\nu}$. Here, $\re{\nu}$ denotes the normal vector on the reference configuration, see Figure \ref{fig:geom_hyp_par_shell} and \ref{fig:mesh_start_hyp_par_shell}.


\begin{figure}
	\centering
	\includegraphics[width=0.238\textwidth]{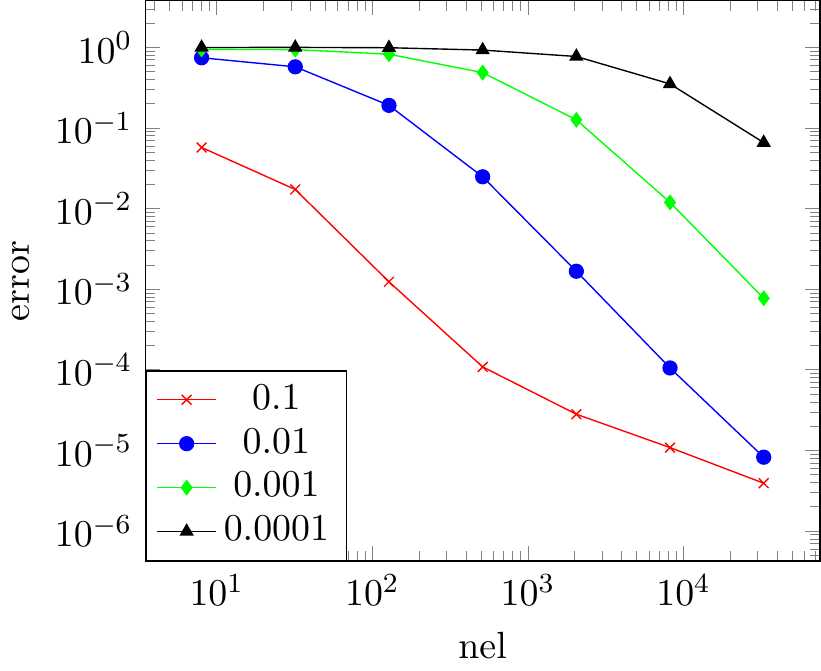}
	\includegraphics[width=0.238\textwidth]{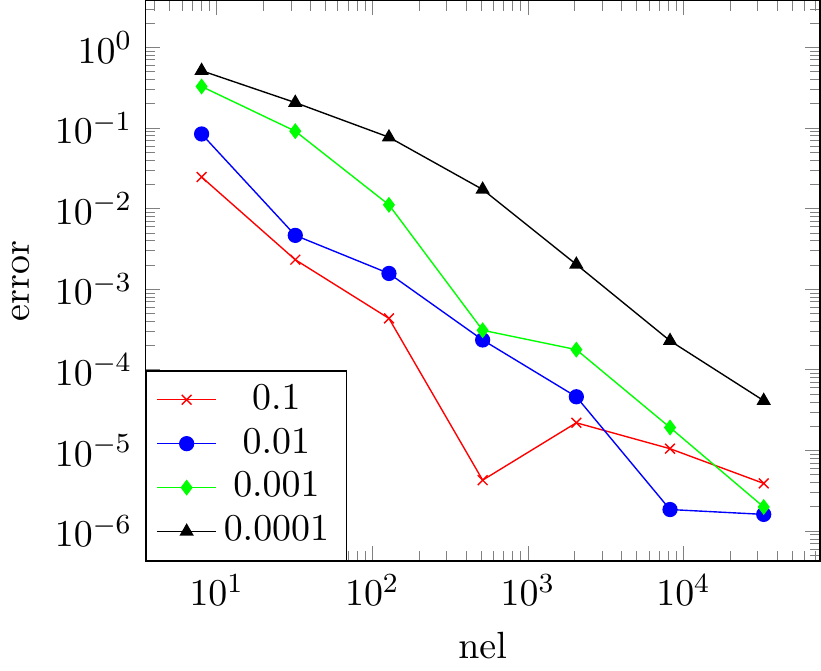}
	\caption{Results for hyperbolic paraboloid, method p2 without and with Regge interpolation.}
	\label{fig:results_hyppar_p2_lin}
\end{figure}

The deflection in z-direction at point $A$ can be seen in Table \ref{tab:results_hyppar_p2_nregge}-\ref{tab:results_hyppar_p2_regge} and the relative error is depicted in Figure \ref{fig:results_hyppar_p2_lin}.

\begin{table}
	\input{Hyperbol_p1_strcTrue_linTrue_regFalse.tab}
	\caption{Results for hyperbolic paraboloid $\times 10^2$, method p2 without Regge interpolation.}
	\label{tab:results_hyppar_p2_nregge}
\end{table}
\begin{table}
	\input{Hyperbol_p1_strcTrue_linTrue_regTrue.tab}
	\caption{Results for hyperbolic paraboloid $\times 10^2$, method p2 with Regge interpolation.}
	\label{tab:results_hyppar_p2_regge}
\end{table}

\subsection{Open hemisphere with clamped ends}
\label{subsec:ne_cyl_clamp}
\begin{figure}
	\centering
	\includegraphics[width=0.22\textwidth]{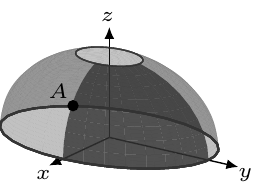}	
	\caption{Geometry for open hemisphere with clamped ends.}
	\label{fig:geom_hemisphere_shell}
\end{figure}

\begin{figure}
	\centering
	\includegraphics[width=0.238\textwidth]{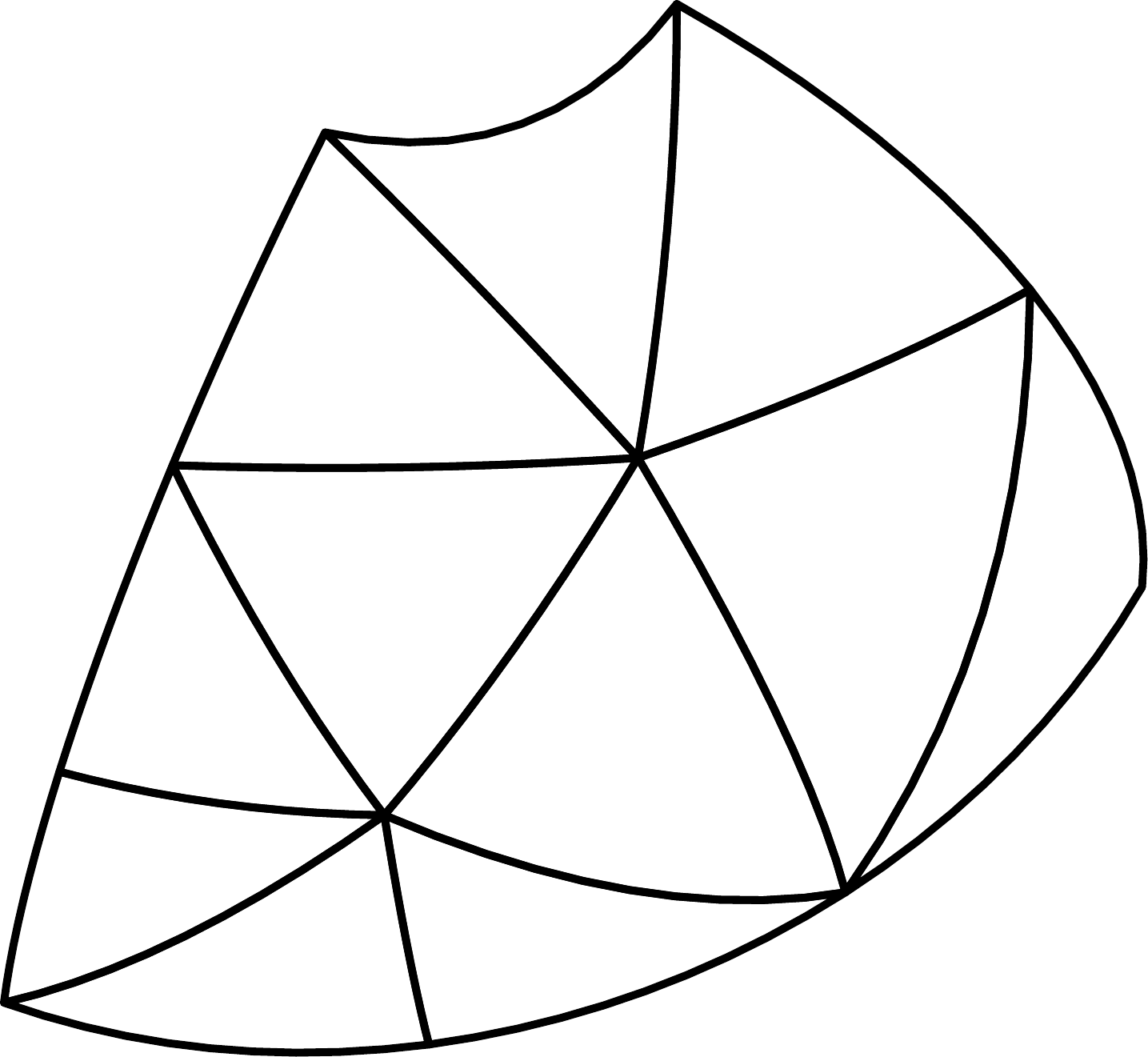}
	\includegraphics[width=0.238\textwidth]{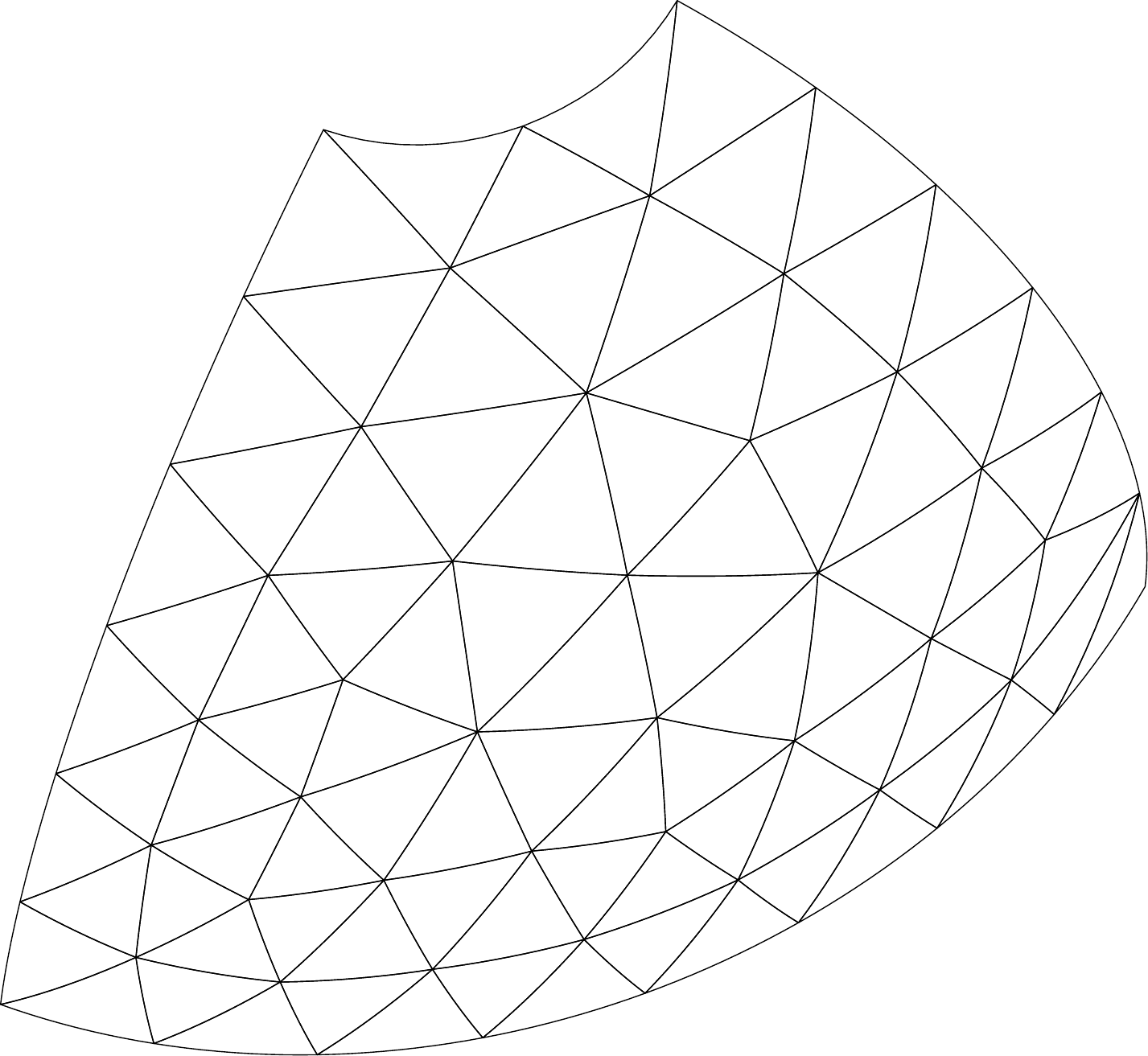}
	\caption{Unstructured meshes with 11 and 84 elements for open hemisphere with clamped ends benchmark.}
	\label{fig:mesh_start_hemisphere_shell}
\end{figure}

As a membrane dominated example an 18$^\circ$ open hemisphere with clamped top and bottom edges is used \cite{Chap11}. Due to symmetry only one fourth of the hemisphere is considered with appropriate symmetry boundary conditions, see Figure \ref{fig:geom_hemisphere_shell} and \ref{fig:mesh_start_hemisphere_shell}. The material and geometric parameters are $R=10$, $E=6.825\times 10^7$, $\nu = 0.3$, $t\in \{0.1, 0.01, 0.001, 0.0001\}$ the volume force density is $P=\frac{t}{10}\cos(2\zeta)\hat{\nu}$, where $\zeta$ denotes the angle between the $x$ and $y$ component.


\begin{figure}
	\centering
	\includegraphics[width=0.238\textwidth]{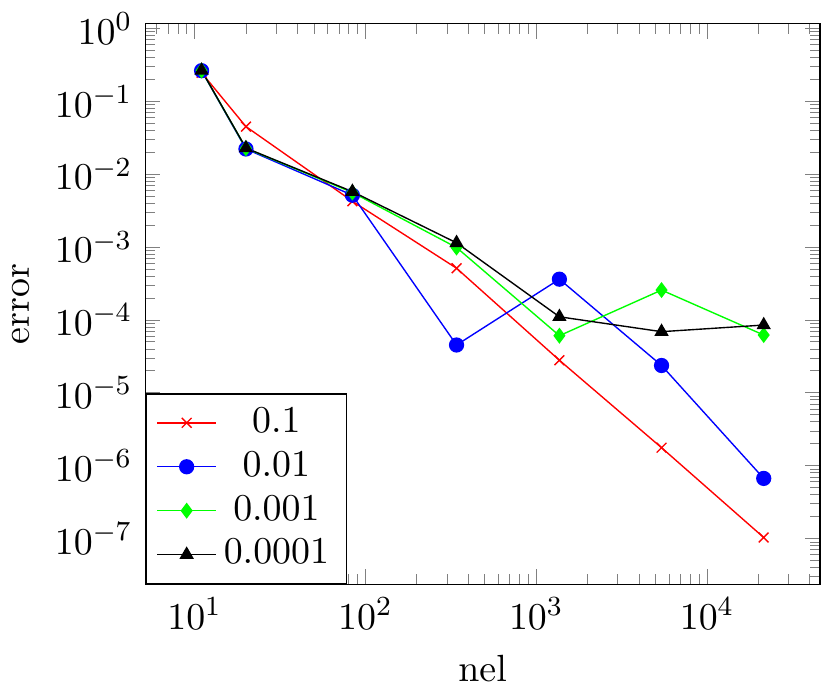}
	\includegraphics[width=0.238\textwidth]{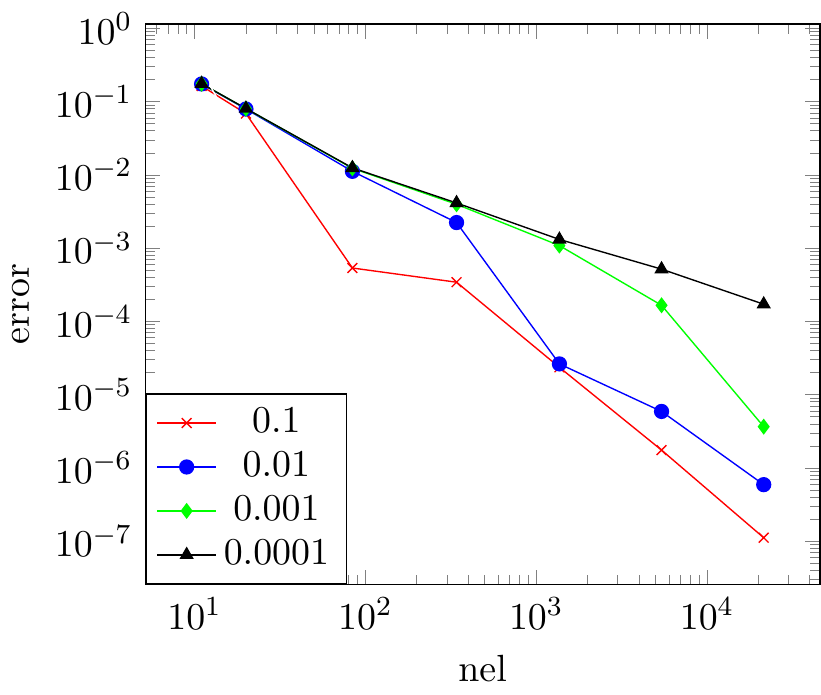}
	\caption{Results for open hemisphere with clamped ends, method p2 without and with Regge interpolation.}
	\label{fig:results_heisphere_p2_lin}
\end{figure}

\begin{table}
	\input{HemisphShell_p1_strcFalse_linTrue_regFalse.tab}
	\caption{Results for open hemisphere with clamped ends $\times 10^5$, method p2 without Regge interpolation.}
	\label{tab:results_hemisphere_p2_nregge}
\end{table}
\begin{table}
	\input{HemisphShell_p1_strcFalse_linTrue_regTrue.tab}
	\caption{Results for open hemisphere with clamped ends $\times 10^5$, method p2 with Regge interpolation.}
	\label{tab:results_hemisphere_p2_regge}
\end{table}

The deflection in $x$-direction at point $A$ is listed in Table \ref{tab:results_hemisphere_p2_nregge}-\ref{tab:results_hemisphere_p2_regge} and Figure \ref{fig:results_heisphere_p2_lin} shows the relative error. As expected the method without Regge interpolation does not lock in the case of inhibited pure bending. Using the interpolation operator yields to only slight deterioration in the convergence rates for smaller thicknesses. However, this effect is marginal compared to the improvements in the bending dominated regime benchmarks. Further, also in this membrane dominated example the interpolation yields to better results for a small amount of elements.

\section{Conclusions}
In this work the Regge interpolation operator was inserted into the membrane energy part relaxing the kernel constraints and avoiding membrane locking for thin shells. For triangular meshes the number of constraints is significantly reduced without deteriorating the membrane stability in the membrane dominated regime. The performance was demonstrated by benchmark examples including membrane and bending dominated cases. A rigorous mathematical proof of uniform convergence independently of the thickness parameter is topic of further research. In contrast to shear locking, which can be already observed for simple plate problems, membrane locking occurs only if curved elements are used and is thus more involved.

\section*{Acknowledgements}
\label{sec:acknowledgements}
The authors acknowledge support from the Austrian Science Fund (FWF) through grant number W 1245.

\bibliographystyle{acm}
\bibliography{cites}
\end{document}

%% file: CylShell_p1_strcFalse_linTrue_regFalse.tab
\begin{tabular}{c|cccc}
 & 0.1 & 0.01 & 0.001 & 0.0001\\
\hline
10 & -2.86758 & -1.46367 & -1.38530 & -1.33837\\
38 & -4.20503 & -1.98840 & -0.74317 & -0.68349\\
160 & -4.57442 & -3.43556 & -1.51590 & -0.37284\\
644 & -4.59672 & -4.36907 & -2.64805 & -0.88707\\
2592 & -4.59775 & -4.46503 & -3.95876 & -2.01521\\
10370 & -4.59781 & -4.47108 & -4.40993 & -3.03108\\
41564 & -4.59781 & -4.47147 & -4.44835 & -4.20861\\
\end{tabular}

%% file: CylShell_p1_strcFalse_linTrue_regTrue.tab
\begin{tabular}{c|cccc}
 & 0.1 & 0.01 & 0.001 & 0.0001\\
\hline
10 & -4.58489 & -4.37366 & -4.21873 & -4.21326\\
38 & -4.59430 & -4.46137 & -4.43641 & -4.42058\\
160 & -4.59764 & -4.46845 & -4.45059 & -4.44846\\
644 & -4.59780 & -4.47137 & -4.44950 & -4.44727\\
2592 & -4.59781 & -4.47150 & -4.45080 & -4.44609\\
10370 & -4.59781 & -4.47150 & -4.45102 & -4.44607\\
41564 & -4.59781 & -4.47150 & -4.45103 & -4.44631\\
\end{tabular}

%% file: CylShell_p1_strcTrue_linTrue_regFalse.tab
\begin{tabular}{c|cccc}
 & 0.1 & 0.01 & 0.001 & 0.0001\\
\hline
8 & -2.29946 & -1.95947 & -1.95357 & -1.95351\\
32 & -3.54279 & -1.85173 & -1.80002 & -1.79858\\
128 & -4.48211 & -2.32588 & -1.78037 & -1.77192\\
512 & -4.58991 & -3.89444 & -1.87506 & -1.76749\\
2048 & -4.59731 & -4.42511 & -2.77129 & -1.78288\\
8192 & -4.59778 & -4.46852 & -4.18701 & -2.00808\\
32768 & -4.59781 & -4.47131 & -4.43271 & -3.36587\\
131072 & -4.59781 & -4.47149 & -4.44987 & -4.33599\\
\end{tabular}

%% file: CylShell_p1_strcTrue_linTrue_regTrue.tab
\begin{tabular}{c|cccc}
 & 0.1 & 0.01 & 0.001 & 0.0001\\
\hline
8 & -4.58733 & -4.48691 & -4.48665 & -4.48665\\
32 & -4.59726 & -4.46548 & -4.45891 & -4.45837\\
128 & -4.59782 & -4.46894 & -4.45131 & -4.45074\\
512 & -4.59781 & -4.47138 & -4.44957 & -4.44759\\
2048 & -4.59781 & -4.47150 & -4.45074 & -4.44619\\
8192 & -4.59781 & -4.47150 & -4.45101 & -4.44602\\
32768 & -4.59781 & -4.47150 & -4.45103 & -4.44630\\
131072 & -4.59781 & -4.47150 & -4.45103 & -4.44633\\
\end{tabular}

%% file: Hyperboloid_p1_strcTrue_linTrue_regFalse.tab
\begin{tabular}{c|cccc}
 & 0.1 & 0.01 & 0.001 & 0.0001\\
\hline
8 & -1.40942 & -0.08393 & -0.00216 & -0.00002\\
32 & -1.99964 & -0.41024 & -0.01227 & -0.00013\\
128 & -2.22356 & -1.23888 & -0.11740 & -0.00182\\
512 & -2.24732 & -1.77470 & -0.61248 & -0.02537\\
968 & -2.24855 & -1.85674 & -1.03445 & -0.07409\\
2048 & -2.24893 & -1.88687 & -1.45795 & -0.20875\\
8192 & -2.24903 & -1.89560 & -1.83745 & -0.91489\\
32768 & -2.24904 & -1.89616 & -1.88895 & -1.64433\\
\end{tabular}

%% file: Hyperboloid_p1_strcTrue_linTrue_regTrue.tab
\begin{tabular}{c|cccc}
 & 0.1 & 0.01 & 0.001 & 0.0001\\
\hline
8 & -2.24057 & -1.92791 & -1.92433 & -1.92429\\
32 & -2.24716 & -1.89856 & -1.89522 & -1.89517\\
128 & -2.24903 & -1.89642 & -1.89304 & -1.89301\\
512 & -2.24904 & -1.89622 & -1.89276 & -1.89271\\
968 & -2.24904 & -1.89621 & -1.89271 & -1.89266\\
2048 & -2.24904 & -1.89620 & -1.89270 & -1.89262\\
8192 & -2.24904 & -1.89620 & -1.89271 & -1.89259\\
32768 & -2.24904 & -1.89620 & -1.89271 & -1.89259\\
\end{tabular}

%% file: UnifBendCyl_p1_strcTrue_linTrue_regFalse.tab
\begin{tabular}{c|cccc}
 & 0.1 & 0.01 & 0.001 & 0.0001\\
\hline
8 & 5.99797 & 5.20307 & 2.43784 & 2.23646\\
16 & 6.00048 & 5.95953 & 4.44977 & 2.89527\\
32 & 6.00005 & 5.99788 & 5.88840 & 4.62870\\
128 & 6.00005 & 5.99970 & 5.96252 & 5.18276\\
384 & 6.00009 & 6.00007 & 5.99700 & 5.80409\\
1024 & 6.00010 & 6.00010 & 5.99987 & 5.98020\\
2560 & 6.00010 & 6.00010 & 6.00008 & 5.99852\\
6144 & 6.00010 & 6.00010 & 6.00010 & 5.99998\\
14336 & 6.00010 & 6.00010 & 6.00010 & 6.00009\\
\end{tabular}

%% file: UnifBendCyl_p1_strcTrue_linTrue_regTrue.tab
\begin{tabular}{c|cccc}
 & 0.1 & 0.01 & 0.001 & 0.0001\\
\hline
8 & 5.99838 & 6.00183 & 6.00186 & 6.00186\\
16 & 5.99925 & 6.00010 & 6.00011 & 6.00011\\
32 & 5.99986 & 6.00007 & 6.00007 & 6.00007\\
128 & 6.00004 & 6.00009 & 6.00009 & 6.00009\\
384 & 6.00009 & 6.00010 & 6.00010 & 6.00010\\
1024 & 6.00010 & 6.00010 & 6.00010 & 6.00010\\
2560 & 6.00010 & 6.00010 & 6.00010 & 6.00010\\
6144 & 6.00010 & 6.00010 & 6.00010 & 6.00010\\
14336 & 6.00010 & 6.00010 & 6.00010 & 6.00010\\
\end{tabular}

%% file: Hyperbol_p1_strcTrue_linTrue_regFalse.tab
\begin{tabular}{c|cccc}
 & 0.1 & 0.01 & 0.001 & 0.0001\\
\hline
8 & 1.72969 & 0.23110 & 0.03778 & 0.00235\\
32 & 1.80316 & 0.38329 & 0.04234 & 0.00063\\
128 & 1.83255 & 0.72759 & 0.12137 & 0.00630\\
512 & 1.83461 & 0.87677 & 0.34994 & 0.04454\\
2048 & 1.83476 & 0.89757 & 0.59515 & 0.13695\\
8192 & 1.83479 & 0.89897 & 0.67297 & 0.38658\\
32768 & 1.83481 & 0.89906 & 0.68057 & 0.55735\\
\end{tabular}

%% file: Hyperbol_p1_strcTrue_linTrue_regTrue.tab
\begin{tabular}{c|cccc}
 & 0.1 & 0.01 & 0.001 & 0.0001\\
\hline
8 & 1.78963 & 0.97481 & 0.90386 & 0.90289\\
32 & 1.83907 & 0.89488 & 0.74326 & 0.71977\\
128 & 1.83561 & 0.89766 & 0.68867 & 0.64244\\
512 & 1.83481 & 0.89928 & 0.68089 & 0.60693\\
2048 & 1.83477 & 0.89911 & 0.68098 & 0.59784\\
8192 & 1.83479 & 0.89907 & 0.68111 & 0.59649\\
32768 & 1.83481 & 0.89907 & 0.68110 & 0.59660\\
\end{tabular}

%% file: HemisphShell_p1_strcFalse_linTrue_regFalse.tab
\begin{tabular}{c|cccc}
 & 0.1 & 0.01 & 0.001 & 0.0001\\
\hline
11 & 1.40550 & 1.41500 & 1.41510 & 1.41510\\
20 & 1.17835 & 1.14513 & 1.14461 & 1.14460\\
84 & 1.13244 & 1.12607 & 1.12551 & 1.12551\\
342 & 1.12707 & 1.12036 & 1.12034 & 1.12034\\
1368 & 1.12762 & 1.11990 & 1.11917 & 1.11919\\
5420 & 1.12765 & 1.12028 & 1.11895 & 1.11899\\
21482 & 1.12765 & 1.12031 & 1.11917 & 1.11897\\
\end{tabular}

%% file: HemisphShell_p1_strcFalse_linTrue_regTrue.tab
\begin{tabular}{c|cccc}
  & 0.1 & 0.01 & 0.001 & 0.0001\\
\hline
11 & 1.31232 & 1.31583 & 1.31587 & 1.31587\\
20 & 1.20520 & 1.20922 & 1.20927 & 1.20927\\
84 & 1.12825 & 1.13285 & 1.13302 & 1.13303\\
342 & 1.12726 & 1.12282 & 1.12369 & 1.12370\\
1368 & 1.12762 & 1.12028 & 1.12045 & 1.12054\\
5420 & 1.12765 & 1.12030 & 1.11942 & 1.11964\\
21482 & 1.12765 & 1.12031 & 1.11923 & 1.11926\\
\end{tabular}